\documentclass{amsart}
\usepackage[active]{srcltx}

\usepackage{amsmath}
 \usepackage{epsfig}
\usepackage{amssymb}
\usepackage{amsfonts}
 \usepackage{latexsym}
 
\newtheorem{Proposition}{Proposition}

  \newtheorem{Lemma}[Proposition]{Lemma}
   \newtheorem{Theorem}[Proposition]{Theorem}
 \newtheorem{Definition}[Proposition]{Definition}

  \makeatletter
\@addtoreset{equation}{section}
\makeatother 

    \def\z{\noindent}  
 \def\Box{{\hfill\hbox{\enspace${\sqre}$}} \smallskip}
    \def\sqr#1#2{{\vcenter{\vbox{\hrule height .#2pt
                             \hbox{\vrule width .#2pt height#1pt \kern#1pt
                                   \vrule width .#2pt}
                             \hrule height .#2pt}}}}
 \def\sqre{\mathchoice\sqr54\sqr54\sqr{4.1}3\sqr{3.5}3}

     \def\CC{\mathbb{C}}
    
    \def\NN{\mathbb{N}}
    
    \def\RR{\mathbb{R}}
    \def\ZZ{\mathbb{Z}}

    \def\bfk{\mathbf{k}}

\begin{document}

\title{Movable singularities of solutions of difference equations in relation
  to solvability, and study of a superstable fixed point}

\author{Ovidiu Costin and Martin Kruskal} \address{
  Department of Mathematics\\Busch Campus-Hill Center\\ Rutgers
  University\\ 110 Frelinghuysen Rd\\Piscataway, NJ 08854\\ e-mail:
  costin\symbol{64}math.rutgers.edu} 
\hyphenation{trans-se-ries trans-as-ymp-tot-ic}

\keywords{Borel summation; Exponential asymptotics; Singularity analysis; Painlev\'e transcendents.\\
 $ ^1$In the sense
stemming from Stokes original papers and the one favored in exponential
asymptotics literature, {\em Stokes} lines are those where a small
exponential is purely real; on an {\em antistokes} line the exponential
becomes purely oscillatory }

\gdef\shorttitle{Movable singularities in nonlinear systems}

\begin{abstract}
 We overview applications exponential asymptotics and analyzable
function theory to difference equations, in defining an analog of the
Painlev\'e property for them and we sketch the conclusions with respect
to the solvability properties of first order autonomous ones.  It turns
out that if the Painlev\'e property is present the equations are
explicitly solvable and in the contrary case, under further assumptions,
the integrals of motion develop singularity barriers. We apply the
method to the logistic map $x_{n+1}=ax_n(1-x_n)$ where it turns out that
the only cases with the Painlev\'e property are $a=-2, 0, 2$ and $4$ for
which explicit solutions indeed exist; in the opposite case an
associated conjugation map develops singularity barriers.
\end{abstract}
\maketitle

\section{Introduction}
\label{Intro}
\z We present an outline of new methods for determining the position and
the type of singularities of a certain kind of solutions of difference
equations \cite{[CK5]} and use this information to perform Painlev\'e
analysis on them. The approach relies on advances in exponential
asymptotics and the theory of analyzable functions
\cite{[Br],[BR2],[C3],[C1], [C2], [CK2],[CK1],[Invent],[CT],
  [CT2],[E3],[E1],[KS],[STL]}. The main concepts are discussed first.

{\em Analyzable functions}. Introduced by Jean
\'Ecalle, these are mostly analytic functions which at singular points
are completely described by {\em transseries}, much in the same way as
analytic functions are represented at regular points by convergent
series.  In contrast with analytic functions (which are not closed under
division) and with meromorphic functions (which fail to be stable under
integration and composition) analyzable ones are conjectured to be
closed under all operations whence the grand picture of this theory,
that all functions of natural origin are analyzable. In particular,
solutions of many classes of differential and difference equations have
been shown to be analyzable.

{\em Transseries}. Also introduced by J. Ecalle, transseries represent
the ``ultimate'' generalization of Taylor series. Transseries are formal
asymptotic combinations of power series, exponentials and logarithms and
contain a wealth of information not only about local but also about
global behavior of functions. One of the simplest nontrivial examples of
a transseries is

\begin{equation}\label{l1}\sum_{k=0}^{\infty}\frac{k!}{x^{k+1}}+Ce^{-x}\ \ \ \ (x\rightarrow +\infty)\end{equation}

\z What distinguishes the expression (\ref{l1}) from a usual asymptotic
expansion is the presence beyond all orders of the principal series of an
exponentially small term which cannot be captured by classical Poincar\'e
asymptoticity.  More general transseries arising in generic ordinary
differential and difference equations are doubly infinite series whose terms
are powers of $x$ multiplied by exponentials (see e.g.  (\ref{transs})). These
transseries can be determined algorithmically and usually diverge factorially,
but can be shown to be Borel summable in a suitable sense in some sector.  In
this sector the function to which the transseries sums has good analyticity
properties and on the edges of the sector typically singularities appear.

{\em The correspondence between the formal universe and actual
functions}. Envisaged by \'Ecalle to be a ``total'' isomorphism
implemented by generalized Borel summation, the correspondence has been
established rigorously in a number of problems including ODEs,
difference equations and PDEs.

\section{{ Difference equations: the isolated movable singularity
    property (IMSP)}} The problem of integrability has been and is a
subject of substantial research. In the context of differential
equations there exist numerous effective criteria of integrability,
among which a crucial role is played by the Painlev\'e test (for
extensive references see e.g. \cite{Kr-Cl}).

An analog of the Painlev\'e property for difference equations turns out
to be more delicate to define.

In the context of solvability, various methods have been proposed by
Joshi \cite{Joshi}, Ablowitz et al. \cite{[AHH]}, Ramani and Grammaticos
(see references in \cite{RG}), Conte and Musette \cite{[CM]}.  See also
\cite{[AHH]} for a comparative discussion of the various approaches in
the literature. None of these is a proper extension of the Painlev\'e
test.  One difficulty resides in continuing the solutions $x(n)$ of a
difference equation, which are defined on a {\em discrete set}, to the
complex plane of the independent variable $n$ in a natural and effective
fashion.  The embedding of $x(n)$ must be done in such a way that {\em
properties are preserved}.  It is important for the effectiveness of the
analysis that this embedding $x(n)${$\sqsubset$}$x(z)$ is natural,
constructive and unique under proper conditions.

It is of course crucial that we are given infinitely many values of the
function $x(n)$; since the accumulation point of $n$ is infinity, the
behavior of $x$ at $\infty$ is key, and then the question boils down to
when infinitely many values determine function, and in which class.

To illustrate a point we start with a rather trivial example. Assume
that $x(n)$ is expressed as a convergent power series in powers of $1/n$

$$
 x_n=\sum_{k=1}^{\infty}c_kn^{-k}  
 $$
 We would then naturally define $x(n)\sqsubset x(z)$ by

$$
 x(z)=\sum_{k=1}^{\infty}c_kz^{-k}
$$

\z for large enough $z$.  Uniqueness is ensured by the analyticity
at $\infty$ of $x(z)$. Still by analyticity, {$\sqsubset$} preserves all
properties.

We cannot rely merely on analytic continuation since there are very few
classes of equations where solutions are given by convergent power
series.  However, as discovered by \'Ecalle and proved in detail very
recently in a general setting by Braaksma \cite{[BR2]}, a wide class of
arbitrary order difference equations admit formal solutions as Borel
summable {\em transseries}.  The class considered by Braaksma is of the
form

\begin{equation}
  \label{DE}
   \mathbf{x}_{n+1}=\mathbf G(\mathbf x_n,n)=\left(\hat\Lambda  +\frac{1}{n}\hat{A}\right)\mathbf{x}_n +\mathbf{g}(n,\mathbf{x}_n)
\end{equation}

\z where $\mathbf{G}$ analytic at $\infty$ in $n$ and at $0$ in $x_n$,
under genericity assumptions \cite{[BR2]}. In particular
a nonresonance condition  is imposed

\begin{equation}\label{condmu}
\mu_m=\bfk\cdot\boldsymbol\mu \mod 2\pi i\end{equation}

\z with $\bfk\in\NN^n$ iff $\bfk=\mathbf{e}_m$.  Transseries solutions
for these equations have many similarities to transseries solutions of
differential equations.  With some $m_1$, $0\le {m_1}\le n$,

    \begin{equation}
      \label{transs}
     \tilde{\bf  x}(n)
=\sum_{{\bf k}\in\NN^{m_1}}{\bf C}^{\bf k} e^{-{\bf k}\cdot
  \boldsymbol{\mu}n}
n^{{\bf k}\cdot {\bf a}}\tilde{\bf y}_\bfk(n)
    \end{equation}
    
    \z In (\ref{transs}) $\boldsymbol{\mu}=(\mu_1,...,\mu_{m_1})$ and
$\mathbf{a}= (a_1,...,a_{m_1})$ depend only on the recurrence (they are
the eigenvalues of $\hat{\Lambda}, \hat{A}$, respectively), ${\bf C}$ is
a free parameter and $\tilde{\bf y}_\bfk$ are formal series in negative
integer powers of $n$, independent on $\mathbf{C}$. The number $m_1$ is
chosen so that all the exponentials in (\ref{transs}) tend to zero in
the chosen sector.
    
    Braaksma showed that $\tilde{\bf y}_\bfk(n)$ are Borel summable
    uniformly in $\bfk$. Let $\mathbf{Y}_\mathbf{k}=\mathcal{L}^{-1}
    \tilde{\bf y}_\bfk(n)$. Then ${\bf Y}_{\mathbf{k}}(p)$ are analytic
    in a neighborhood of $\RR^+$ (in fact in a larger sector). Defining

\begin{equation}\label{-2}
{\bf y}_{\mathbf{k}}=\int_0^{\infty}e^{-np}{\bf
    Y}_{\mathbf{k}}(p)dp\end{equation}

\z  we have uniform estimates $|\mathbf y_\mathbf
    k|<A^\mathbf k$ and thus the series
\begin{equation}
  \label{eq:transsS1}
  {\bf x}(n)
=\sum_{{\bf k}\in\NN^{m_1}}{\bf C}^{\bf k} e^{-{\bf k}\cdot
  \boldsymbol{\mu} n}
n^{{\bf k}\cdot {\bf a}}{\bf y}_\bfk(n)
\end{equation}
\z is classically convergent for large enough $n$. Braaksma showed that
${\bf x}(n)$ is an actual solution of (\ref{DE}).  It is natural to
replace $n$ by $z$ in (\ref{-2}) and define:
\begin{equation} \label{eq:transsS}
  {\bf x}(z)
=\sum_{{\bf k}\in\NN^{m_1}}{\bf C}^{\bf k} e^{-{\bf k}\cdot
  \boldsymbol{\mu} z}
z^{{\bf k}\cdot {\bf a}}{\bf y}_\bfk(z)
\end{equation}
\z If $z$ and all constants are real and $\mu_i<0$, the functions
(\ref{eq:transsS}) are special cases of \'Ecalle's {\em analyzable}
functions. As explained before we are allowing for $z$ and constants to
be complex, under restriction $\Re({\bf k}\cdot \boldsymbol{\mu} z)>0$.

It is crucial that the values of $\bf y (n)$ for all large enough $n$ uniquely
determine the expansion.  In \cite{[CK5]} it is shown that under suitable
conditions, {\em two distinct analyzable functions cannot agree on a set of
  points accumulating at infinity.} Below is a simplified version of a theorem
in \cite{[CK5]}.

\begin{Theorem}[\cite{[CK5]}]  {Assume 
    
$$
  \label{eq:nonres}
  \left(\mathbf{Z}\cdot \boldsymbol\mu=0 \mod 2\pi i \ \ \mbox{with}\ \ \mathbf{Z}\in\ZZ^p 
  \right)
\Leftrightarrow \mathbf{Z}=0
$$
and $ {\bf x}(z)=\sum_{{\bf k}\in\NN^{m_1}}{\bf C}^{\bf k} e^{-{\bf
    k}\cdot \boldsymbol{\mu} z} z^{{\bf k}\cdot {\bf a}}{\bf y}_\bfk(z)$.
If $ {\bf x}(n)=0$ for all large enough $n\in\NN$, then $
{\bf x}(z)$ is identically zero.}

\end{Theorem}

Analyzable functions behave in most respects as analytic functions. Among the
common properties, particularly important is the uniqueness of the extension
from sets with accumulation points, implying the principle of permanence of
relations.  Under the assumptions \cite{[BR2]} that ensure Borel summability,
we have the following.

\begin{Definition}
  \label{D3}  If 

$${\bf x}(n)
=\sum_{{\bf k}\in\NN^{m_1}}{\bf C}^{\bf k} e^{-{\bf k}\cdot
  \boldsymbol{\mu} n}
n^{{\bf k}\cdot {\bf a}}{\bf y}_\bfk(n)$$

\z  then we call 

$${\bf x}(z)
=\sum_{{\bf k}\in\NN^{m_1}}{\bf C}^{\bf k} e^{-{\bf k}\cdot
  \boldsymbol{\mu} z}
z^{{\bf k}\cdot {\bf a}}{\bf y}_\bfk(z)$$

\z   the { analyzable embedding} of $\mathbf{x}(n)$
  from $\NN$ to a sector in $\CC$.

\end{Definition}

\z Having now a suitable procedure of analytic continuation, we can
define the isolated movable singularity property, an extension of the
Painlev\'e property to difference equations, in a natural fashion.  In
analogy with the case of differential equations we require that all
solutions are free from ``bad'' movable singularities.

\begin{Definition}
  A
  difference equation has the IMSP if all {\em movable}
  singularities of all its solutions are {\bf isolated.}
\end{Definition}
{\bf Notes.} (i) We use the common convention that isolated singularity
exclude branch points, clusters of poles and barriers of singularities.

(ii) To determine the singularities,  transasymptotic matching
methods introduced for differential equations in \cite{[Invent]}, can be
extended with little changes to difference equations.

\section{Classification of some difference
    equations with IMSP. Solvability.}

We look at autonomous difference equations of the form $x_{n+1}=G(x_n)$
where $G$ is meromorphic and has attracting fixed points. A more general
analysis is given in \cite{[CK5]}.  We then write
\begin{equation}
  \label{eq:01}
    x_{n+1}=G(x_n):=ax_n+F(x_n)
\end{equation}

\z and restrict for simplicity to the case $F(0)=F'(0)=0$ and
$0<|a|<1$. There is a one-parameter family of solutions presented as
simple transseries {\em convergent} for large enough $n$, of the form

\begin{equation}
  \label{eq:2}
 x_n=x_n(C)=\sum_{k=1}^\infty e^{nk\ln a} C^k D_k
\end{equation}

\z for given values of $D_k$, independent of $C$. The analyzable
embedding of $x$, cf. Definition \ref{D3}, reads
\begin{equation}
  \label{eq:2z} x(z)=x(z;C)=\sum_{k=1}^\infty e^{z k\ln a} C^k D_k
\end{equation}

\z which is analytic for large enough $z$. To look for the IMSP, we find the
properties of $x(z)$ beyond the domain of convergence of (\ref{eq:2z}), and
find the singular points of $x(z)$.

{\bf Note.} Because equation (\ref{eq:01}) is nonlinear, although
(\ref{eq:2z}) has one continuous parameter, there may be more solutions.
This issue is addressed in (\cite{[CK5]}).

\subsection{Embedding versus properties of the conjugation map} By the
Poincar\'e conjugation theorem applied to $x_{n+1}=G(x_n)$ there exists
a unique map $\phi$ with the properties

\begin{equation}
  \label{assum1}
  \phi(0)=0,\ \  \phi'(0)=1\ \ \mbox{and} \ \phi \ \mbox{analytic at }\  0
\end{equation}

\z and such that $x_n=\phi(X_n)$ implies $X_{n+1}=aX_n$.  The map $\phi$
is a conjugation map of $x_{n+1}$ with its linearization $X_{n+1}$. We
have $X_n=a^n X_0=Ca^n$.

We obtain an extension of $x$ from $\NN$ to $\CC$ through

\begin{equation}
  \label{eq:con1}
  x_n=\phi(Ca^n)\ \sqsubset\  x(z):=\phi(Ca^z)
\end{equation}

Then the conjugation map satisfies the equation
$$
  \label{eq:conjug}
  \phi(az)=G(z)=a\phi(z)+F(\phi(z))
$$

As expected by the uniqueness of the embedding in Definition~\ref{D3} we
have the following.

\begin{Lemma}\label{R1}
   (i) For equations of the type (\ref{eq:01}) under the given
    assumptions, the embeddings through analyzability and via the
    conjugation map, cf. (\ref{eq:con1}) agree.
    
    (ii) The solutions $x(z;C)$ have only isolated movable singularities
    iff $\phi$ has only isolated singularities.
\end{Lemma}

\z {\bf Corollary} {\em A necessary condition for (\ref{eq:01}) to have the IMSP
is that the conjugation map $\phi$, around every stable fixed point of
$G$ and of those of its iterates which are meromorphic (of the form
$G^{[m]}=G\circ\cdots\circ G$ $ m$ times) extends analytically into the
complex plane, except for isolated singularities.}

\begin{figure}
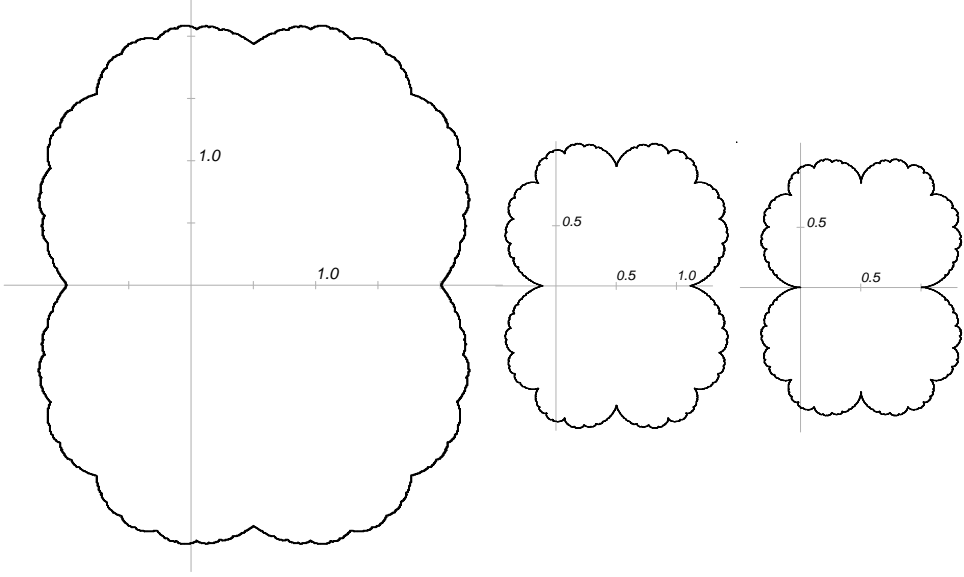
{\label{Fig.1}}
\begin{picture}(500,250)\hskip -3cm
\epsfig{file=050.ps, height=10.10cm}
\end{picture}
\vskip  -2.1cm 
\begin{picture}(100,000)
\epsfig{file=090.ps, height=5cm}
\end{picture}
\vskip -0.4cm \hskip 6.5cm
\begin{picture}(100,000)
\epsfig{file=100.ps, height=5cm}
\end{picture}
\vskip 3cm

\caption{Julia sets  for $G=ax(1-x)$ for $a=0.5$, $a=0.9$ and $a=1$ respectively. The set $\mathcal{K}_p$ is the interior of these curves (graphs generated using C and MapleV)}
\end{figure}

\subsection{Autonomous equations with the IMSP}

In (\cite{[CK5]}) we classify the equations of the form (\ref{eq:01})
with respect to the IMSP.  In a way analog to the case of ODEs of first
order, only Riccati equations have this property. On the other hand,
these equations are explicitly solvable, and thus autonomous equations
of the first order which have the IMSP can be solved in closed form.

\begin{Theorem}\label{T1}
  The equation (\ref{eq:01}) under the given assumptions {\em fails} to
  have the IMSP {\em unless}, for some $c\in\CC$,
\begin{equation}
\label{(**)}
G(z)=\frac{az}{1+cz}
\end{equation}
\z In case (\ref{(**)} we have: $1/x_n=a^{-n}(C-c/(a-1))+c/(a-1)$
\end{Theorem}

It is also very interesting to see that in the case of failure of the
IMSP the analytic properties of the solutions preclude the existence of
nice constants of motion. This is discussed in the next section.

\subsection{Case of failure of IMSP} We illustrate this situation when 
$G$ polynomial.  The surprising conclusion is that for polynomial $G$
without the IMSP, constants of motion (defined as functions $C(x,n)$
constant along trajectories) develop barriers of singularities along
some fractal closed curves $\partial \mathcal{K}_p$, see below. In a
neighborhood of the origin, the function $\phi$ is invertible, and thus
for small $x$ we derive from $x_n=\phi(a^nC)$ that $ C=a^{-n}Q(x_n)$
where $Q=\phi^{-1}$.

 In Fig. 1. we depict the Julia set $\partial \mathcal{K}_p$ of a simple
map. In that case, in the compact set bounded by $\mathcal{K}_p$
consists in the initial conditions for which the solution of the
iteration converges to zero.  The Julia set is a closed curve of
nontrivial fractal dimension.  For a comprehensive discussion of Julia
sets and iterations of rational maps see \cite{Beardon}.

\begin{Theorem}\label{Barrier}
 Assume $G$ is a polynomial map with an attracting fixed point at the
  origin.  Denote by $\mathcal{K}_p$ (cf. Fig. 1) the maximal connected
  component of the origin in the Fatou set of $G$.

 Then the domain of analyticity of $Q$ is $\mathcal{K}_p$, and
  $\partial \mathcal{K}_p$ is a barrier of singularities of $Q$.
\end{Theorem}

The proofs of these results can be found in \cite{[CK5]}. The logistic
map discussed in relative detail in the next section represents a very
simple illustration of some of the relevant phenomena.

\section{Analysis of the logistic map at the superstable fixed
  point infinity}

 We show that the equation $x_{n+1}=ax_n(1-x_n)$ has the IMSP iff
$a\in\{-2,0,2,4\}$. The case $a=0$ needs no analysis. Otherwise, taking
$y=1/x$ we get

\begin{equation}
  \label{eq:111}
 y_{n+1}=\frac{y_{n}^2}{a(y_n-1)}
\end{equation}
  
\z For small $y_0$, the leading order form of equation (\ref{eq:111}) is
$y_{n+1}=-a^{-1}{y_{n}^2}$ whose solution is $-y_0^{2^n}a^{-2^{n-1}-1}$. It is
then convenient to seek solutions of (\ref{eq:111}) in the form
$y_n=F(y_0^{2^n}a^{-2^{n}})$ whence the initial condition implies $F(0)=0$,
$F'(0)=-a$. Denoting $y_0^{2^n}a^{-2^{n}}=z$, the functional relation
satisfied by $F$ is

\begin{equation}
  \label{eq:112}
 F(z^2)=\frac{F(z)^2}{a(F(z)-1)} \ \ \ \ (F(0)=0,\ F'(0)=-a)
\end{equation}

\begin{Lemma}\label{L11}
(i) There exists a unique analytic function $F$ in the neighborhood of the
origin satisfying (\ref{eq:112}) and such that $F(0)=0$ and
$F'(0)=-a$. This $F$ has only isolated singularities in $\CC$ if and
only if $a\in\{-2,2,4\}$. In the latter case, the equations
(\ref{eq:112}) and (\ref{eq:111}) can be solved explicitly (see
(\ref{114.5}), (\ref{119}) and (\ref{120})). 

If $a\not\in\{-2,2,4\}$ then the unit disk is a barrier of singularities
of $F$.

(ii) If $a\not\in\{-2,2,4\}$ then $\partial \mathcal{K}_p$ is a barrier of
singularities of the constant of motion $Q=F^{-1}$ (cf.
Theorem~\ref{Barrier}).

\end{Lemma}

\subsection{Proof of Lemma \ref{L11}} 

\subsubsection{Analyticity at zero} We write 

\begin{equation}
  \label{eq:113}
  F(z)=\frac{a}{2}F(z^2)-\frac{1}{2}\sqrt{a^2 F(z^2)^2-4aF(z^2)}
\end{equation}

\z (with the choice of branch consistent with $F(0)=0, F'(0)=-a$), and take
$F(z)=-az+h(z)$. This leads to the equation for $h$

\begin{equation}
  \label{eq:1131}
  h(z)=az+\frac{a}{2}(h(z^2)-az^2)-\frac{1}{2}\sqrt
{4a^2z^2-4ah(z^2)+a^2(h(z^2)-az^2)^2}=\mathcal{N}(h)
\end{equation}

\z It is straightforward to show that, for small $\epsilon$, $\mathcal{N}$
is contractive in the ball of radius $|a|^2$ in the space of functions
of the form $h(z)=z^2u(z)$, where $u$ is analytic in the disk
$D=\{z:\|z\|<\epsilon\}$ with the norm $\|h\|=\sup_{D}|u|$.  The corresponding
$F$ is analytic for small $z$ and is a conjugation map between
(\ref{eq:111}) and its small-$y$ approximation. 

\subsubsection{} Now the question is to determine the singularities of $F$ in the complex
plane. Inside the unit disk we can inductively continue $F$ analytically
through (\ref{eq:113}) as follows. If $F$ is analytic in a disk of radius
$r<1$ then (\ref{eq:113}) provides the analytic continuation in the disk of
radius $\sqrt{r}$ if $F(z^2)$ is not zero or $4/a$. Because $F$ is assumed to
be analytic in the disk of radius $r$, $F(z^2)$ cannot vanish; otherwise by
(\ref{eq:112}) it would vanish infinitely often in a neighborhood of zero,
which is inconsistent with its analyticity at $z=0$ and with $F'(0)=-a$.  But
the equation $F(z^2)=4/a$ does have solutions in the unit disk if $|a|>5$.
Indeed, with $|a|>5$ it is immediate to show $F(z^2)=4/a$ implies
$F(z^{2^n}):=\epsilon_n\rightarrow 0$ as $n\rightarrow\infty$. Since $F$ is,
by the condition $F'(0)=-a$, bijective near $z=0$, the equation
$F(z)=\epsilon_n$ has a (unique) solution, $z=z_0$ for sufficiently large $n$.
Now, if $z_1$ is such that $z_1^{2^n}=z_0$ we have $F(z_0)=4/a$. But then $F$
has a square root branch point at ${z_1}$ since $F'(z)$ cannot vanish anywhere
inside the unit disk, otherwise again by (\ref{eq:112}), $F'(z)$ would vanish
infinitely often near the origin.

Thus for $|a|>5$ there is no $F$ with only isolated singularities. In fact, as
is seen in the last paragraph of the proof, in this case the unit disk is a
barrier of singularities of $F$.

We now consider the case where $|a|<5$.

\subsubsection{} We claim that unless $a=-2$, the point $z=1$ is a singular point of $F$.
Indeed, a Taylor series expansion $F=\sum_{k=0}^\infty c_k(z-1)^k$ gives
$c_0=0$ or $c_0=a/(a-1)$. It is straightforward to see that $c_0=0$ implies
$c_k=0$ for all $k$ which contradicts $F'(0)=-a$. 

Therefore $c_0=a/(a-1)$ in which case direct calculation shows that, unless
$a=-2(2^k-1)$ for some $k\in\ZZ$, all $c_k$ for $k\ge 1$ are zero which is not
possible since $F(0)=0$. If $k>1$ then $|a|>5$.

Therefore, $k=1$ whence $a=-2$. With this value of $a$, (\ref{eq:112}) has the
explicit solution

\begin{equation}\label{114.5} F(z)=\frac{2z}{z^2+z+1}\end{equation}
It remains to look at the cases when $z=1$ is a singular point of
$F$.

\subsubsection{} If  $z=1$ is a meromorphic point of  $F$, then we obtain from (\ref{eq:112}):  

\begin{equation}\label{116}
\lim_{z\rightarrow 1}\frac{aF(z^2)}{F(z)}=1\end{equation}

\z If $F(z)=\sum_{k=-p}^{\infty}c_k(z-1)^k$ is the Laurent series of $F$ at
$z=1$, then (\ref{116}) implies

\begin{equation}
  \label{eq:117}
  a=2^p
\end{equation}

\z and, since $|a|<5$, we must have $a\in\{2,4\}$. 

\z For $a=2$ (\ref{eq:112}) has the explicit solution

\begin{equation}\label{119}F(z)=\frac{2z}{z-1}\end{equation}

\z For $a=4$ the solution of (\ref{eq:112}) is

\begin{equation}\label{120}F(z)=-\frac{4z}{(z-1)^2}\end{equation}

\subsubsection{}  If $F$ is not meromorphic at $z=1$ then $F$ cannot be continued
analytically beyond the unit circle:

 Assume first that $F$ is analytic in the open unit disk $D_1$. By
(\ref{eq:112}) we have

\begin{equation}\label{122}F(z^{2^n})=R_n(F(z))\end{equation}

\z where $R_n$ is a rational function. Thus if $z_0^{2^n}=1$ and 
$F$ is analytic at $z_0$ then it is meromorphic at $z=1$, while the set
of points such that $z_0^{2^n}=1$ is dense on the unit circle.

The last case to analyze is $R<1$, where $R$ is the (maximal) radius of
analyticity of $F$.  Equation (\ref{eq:113}) shows that on the disk of radius
$R$ the singularities of $F$ are square root branch points. On the other hand,
if $z_0\in D_1$ is an algebraic branch point and $z_1^2=z_0$ then $z_1$ is
also an algebraic branch point as can be seen by continuing (\ref{eq:112})
around a small circle near the branch point $z^2=z_0$.  It is then easy to see
that if $R<1$ the branch points accumulate densely towards the unit circle.
$\Box$

(ii) The proof follows the lines of \ref{[Ck5]} so we only sketch the
main steps. We note that by (i) $F$ is invertible near the origin,
so $Q$ is analytic near the origin, where it satisfies the relation
$Q^2(z)=Q(z^2/(az-1)$.


\begin{thebibliography}{99}



\bibitem{[AHH]} M J Ablowitz, R Halburd, and B Herbst {\em On the
extension of the Painlev\'e property to difference equations}
Nonlinearity {\bf 13} pp. 889--905 (2000).



\bibitem{Arnold} D V Anosov and V I Arnold eds. Dynamical Systems I
Springer-Verlag  (1988).

\bibitem{Baker} I N  Baker {\em Fixpoints of polynomials and rational
  functions} {J. London Math. Soc.}  {\bf 39} pp. 615--622 (1964).

\bibitem{Balser} W Balser {\em From divergent power series to
    analytic functions, Springer-Verlag, 1994.}

\bibitem{Beardon} A F Beardon {\em Iteration of Rational Functions}
Springer Verlag, New York (1991).

\bibitem{benderorszag} C Bender and S Orszag, {\em Advanced
Mathematical Methods for scientists and engineers}, McGraw-Hill, 1978.

\bibitem{[Be1]} M V        Berry, \emph{Uniform asymptotic smoothing of Stokes's discontinuities}, Proc. R. Soc. Lond. A 422, 7-21 (1989).

\bibitem{[BH]} M V        Berry, C J Howls \emph{Hyperasymptotics}, Proc. R. Soc. Lond. A 430, 653-668 (1990).




\bibitem{BorelM} E Borel {\em Le\c cons sur les fonctions monog\`enes,
Gauthier-Villars}, Paris (1917).

\bibitem{[Br]} B L J     Braaksma, \emph{Multisummability of formal power series solutions of nonlinear meromorphic differential equations}, Ann.  Inst. Fourier, Grenoble,{\bf 42}, 3, 517-540 (1992).

\bibitem{[BR2]} B L J  Braaksma, {\em Transseries for a class of
nonlinear difference equations} (To appear in Journ. of Difference
Equations and Applications). 

  
\bibitem{[CM]} R Conte and M Musette {\em Rules of discretization for
Painlev\'e equations} Theory of Nonlinear Special Functions (Montreal
13--17 May 1996) ed. L Vinet and P Winternitz (Berlin:Springer).

 
\bibitem{[C1]} O           Costin, \emph{Exponential asymptotics, transseries, and generalized Borel summation for analytic, nonlinear, rank-one systems of ordinary differential equations}, Internat. Math. Res. Notices
    no. 8, 377--417 (1995).$^1$
  



    
  \bibitem{[C2]} O  Costin, \emph{On Borel summation and Stokes
      phenomena for rank one nonlinear systems of ODEs}, Duke Math. J.
    Vol. 93, No.2 (1998).$^1$


  
\bibitem{[CK2]} O  Costin and M D Kruskal, \emph{On optimal truncation
    of divergent series solutions of nonlinear differential systems;
    Berry smoothing}, Proc. R. Soc. Lond. A {\bf 455}, 1931-1956
  (1999).$^1$

\bibitem{[CK1]} O          Costin and M D Kruskal, \emph{Optimal uniform estimates and rigorous asymptotics beyond all orders for a class of ordinary differential equations},  Proc. Roy. Soc. London Ser. A 452, no.
    1948, 1057--1085 (1996).\setcounter{footnote}{0}\footnote{Available online at http://www.math.rutgers.edu/\raisebox{-0.5em}{\large\symbol{126}}costin}



\bibitem{[CK5]} O          Costin and M D Kruskal, {\em Equivalent of the Painlev\'e property 
for certain classes of difference equations and study of their
solvability} (preprint 2001).$^1$.



\bibitem{[Invent]} O  Costin and R D Costin, {\em On the location and
type of singularities of nonlinear differential systems} (Inventiones
Mathematicae \textbf{145}, 3, pp 425-485 (2001))).$^1$

\bibitem{[CC2]} O          Costin and R D Costin, \emph{Singular normal form for the Painlev\'e equation P$_{\rm I}$}, Nonlinearity Vol 11, No. 5 pp.  1195-1208 (1998).$^1$


\bibitem{[CT]} O           Costin and S Tanveer, {\em Existence and uniqueness for a class of nonlinear higher-order partial differential equations in the complex plane} (CPAM Vol. LIII, 1092---1117 (2000)).$^1$

  
\bibitem{[CT2]} O  Costin and S Tanveer, {\em Existence and uniqueness
of solutions of nonlinear evolution systems of n-th order partial
differential equations in the complex plane} (submitted).$^1$



\bibitem{[CKbook]} O Costin and M D Kruskal, {\em Analyzable functions,
transseries and applications}, book in preparation for Springer-Verlag.


\bibitem{[C3]} O           Costin, \emph{Correlation between pole location and asymptotic behavior for Painlev\'e I solutions}, Comm.  Pure and Appl. Math. Vol. LII, 0461-0478 (1999).




\bibitem{Deift} P A Deift and X Zhou {\em A steepest descent method for
    oscillatory Riemann-Hilbert problems. Asymptotics for the MKdV
    equation}
Ann. of Math. (2)  137 no 2 pp. 295--368 (1993).



\bibitem{[E3]} J  \'Ecalle, {\em Finitude des cycles limites et
acc\'el\'ero-sommation de l'application de retour, Preprint 90-36 of
Universite de Paris-Sud (1990).}
\bibitem{[E1]} J         \'Ecalle, {\em Fonctions Resurgentes, Publications Mathematiques d'Orsay (1981).}

\bibitem{[E2]} J         \'Ecalle, {\em in Bifurcations and periodic
    orbits of vector fields NATO ASI Series, Vol. 408 (1993).}

\bibitem{FlN} H Flashka and A C Newell {\em Monodromy and spectrum
    preserving transformations}
Commun. Math. Phys. {\bf 76} pp. 65--116 (1980).

\bibitem{Fuchs} L Fuchs {\em Sur quelques \'equations diff\'erentielles 
lin\'eaires du second ordre} C. R. Acad. Sci., Paris {\bf 141}
pp. 555-558 (1905).

\bibitem{Gonzales} M O Gonz\'ales {\em Complex Analysis: Selected
    Topics}
Marcel Dekker Inc., NY, Basel, Hong Kong (1991).

\bibitem{Koval} S Kowalevski {\em Sur le probl\`eme de la rotation d'un
corps solide autour d'un point fixe}, Acta Math., {\bf 12} H.2,
pp. 177--232 (1889)

\bibitem{Koval2} S Kowalevski {\em M\'emoire sur un cas particulier du
probl\`eme d'un corps pesant autour d'un point fixe, o\`u l'integration
s'effectue \`a l'aide de fonctions ultraelliptiques du temps} M\'emoires
pr\'esent\'es par divers savants \`a l'Academie des Sciences de
l'Institut National de France, Paris 31 pp. 1--62 (1890)

117--232.
\bibitem{Kr-Cl} M D Kruskal and P A Clarkson {\em The
    Painlev\'e-Kowalevski and poly-Painlev\'e tests for integrability}
    Stud. Appl. Math. 86 no. 2, pp.  87--165 (1992).


  
\bibitem{[KS]} M D Kruskal and H Segur, {\em Asymptotics beyond all
    orders in a model of crystal growth} , Studies in Applied
  Mathematics 85(2), 129-181 (1991)

\bibitem{Gambier} B Gambier {\em Sur les \'equations diff\'erentielles du
  second ordre et du premier degr\'e dont l'int\'egrale g\'en\'erale
est \`a points critiques fixes} Acta Math. {\bf 33} pp 1--55 (1910).

\bibitem{RG} B Grammaticos and A Ramani {\em Discrete Painlev\'e
equations: derivation and properties} In ``Applications of Analytic and
Geometric Methods to Nonlinear Differential Equations'', ed. P A
Clarkson, NATO ASI Series C, pp 299--313 (1993).

\bibitem{Its} A R Its, A S Fokas and A A Kapaev {\em On the asymptotic
  analysis of the Painlev\'e equations via the isomonodromy method}
Nonlinearity 7
no. 5, pp. 1921--1325 (1994).


\bibitem{Joshi} N Joshi {\em Irregular singular behaviour in the first
    discrete Painlev\'e equation}, in Symmetries and Integrability of
  Difference Equations III, D Levi and O Ragnisco (eds), pp. 237--243,
  CRM Proc.  Lecture Notes, 25, Amer. Math. Soc., Providence, RI,
  (2000).

\bibitem{Levy} H Levy and F Lessman {\em Finite Difference Equations}, Dover
  Publications Inc. New York, 1992.
  
\bibitem{Nevanlinna} R Nevanlinna  {\em Le th\'eor\`eme de Picard-Borel et
  la th\'eorie des fonctions m\'eromorphes}, Chelsea Pub. Co., New York,
  (1974).  
  
\bibitem{Painleve1} P Painlev\'e {\em M\'emoire sur les \'equations
    diff\'erentielles dont l'integrale g\'en\'erale est uniforme} Bull.
  Soc. Math. France {\bf 28} pp 201--261 (1900).

  
\bibitem{Painleve2} P Painlev\'e {\em Sur les \'equations
    diff\'erentielles du second ordre et d'ordre sup\'erieur dont
    l'integrale g\'en\'erale est uniforme} Acta Math. {\bf 25} pp. 1--85
  (1902). 

\bibitem{Rudin} W  Rudin {\em Real and Complex
  Analysis}, McGraw-Hill (1987).


\bibitem{Steinmetz} N Steinmetz {\em Rational iteration. Complex Analytic
 Dynamical Systems} Walter de Gruyter, Berlin; New York (1993).


\bibitem{[STL]} H Segur, S  Tanveer and H Levine, ed. {\em
    Asymptotics Beyond all Orders}, Plenum Press (1991).


\bibitem{[Va]} V S       Varadarajan, \emph{Linear meromorphic differential equations: a modern point of view}. Bull. Amer. Math. Soc. (N.S.)  33, no. 1, 1--42 (1996).

\bibitem{Wasow} W Wasow, {\em Asymptotic expansions for ordinary
differential equations}, Interscience Publishers, 1968.





\end{thebibliography}
\end{document}